\documentclass[MSNbibl,number,citesort,seceqn,dvips]{arxbj}
\usepackage{upgreek}
\usepackage{mathrsfs}

\aid{0}
\volume{19}
\issue{2}
\pubyear{2013}
\firstpage{610}
\lastpage{632}
\doi{10.3150/11-BEJ406}

\makeatletter
\def\tsfrac#1#2{\frac{#1}{#2}}

\newcommand{\IZ}{\mathbb{Z}}
\newcommand{\IR}{\mathbb{R}}

\def\dtv{{d_{\mathrm{TV}}}}
\def\dloc{{d_{\mathrm{loc}}}}

\def\dk{{d_{\mathrm{K}}}}

\newcommand{\law}{\mathscr{L}}
\newcommand{\eqlaw}{\stackrel{\mathscr{D}}{=}}

\newcommand{\eqref}[1]{(\ref{#1})}
\newproclaim{definition}{Definition}[section]
\newremark{remark}[definition]{Remark}
\newtheorem{theorem}{Theorem}[section]
\newtheorem{proposition}[theorem]{Proposition}
\newtheorem{lemma}[theorem]{Lemma}

\newcommand{\IE}{\mathbb{E}}
\newcommand{\IP}{\mathbb{P}}
\newcommand{\eq}{\eqref}
\newcommand{\Var}{\operatorname{Var}}
\newcommand{\Exp}{\operatorname{Exp}}
\newcommand{\Ge}{\operatorname{Ge}}
\newcommand{\Be}{\operatorname{}\operatorname{Be}}
\newcommand{\eps}{\varepsilon}
\newcommand{\D}{\Delta}
\newcommand{\I}{\mathrm{I}}
\newcommand{\ahalf}{{\frac{1}{2}}}

\def\mid{\vert}
\def\bbbmid{\vert}

\makeatother

\begin{document}
\begin{frontmatter}

\title{Total variation error bounds for geometric approximation}
\runtitle{Total variation error bounds for geometric approximation}

\begin{aug}
\author[1]{\fnms{Erol A.} \snm{Pek\"oz}\thanksref{1}\ead[label=e1]{pekoz@bu.edu}},
\author[2]{\fnms{Adrian} \snm{R\"ollin}\thanksref{2}\ead[label=e2]{adrian.roellin@nus.edu.sg}}
\and
\author[3]{\fnms{Nathan} \snm{Ross}\corref{}\thanksref{3}\ead[label=e3]{ross@stat.berkeley.edu}}
\runauthor{E.A. Pek\"oz, A. R\"ollin and N. Ross} 
\address[1]{School of Management, Boston University, 595 Commonwealth Avenue
Boston, MA~02215, USA. \printead{e1}}
\address[2]{Department of Statistics and Applied Probability, National
University of Singapore, 6 Science Drive 2, Singapore 117546,
Singapore. \printead{e2}}
\address[3]{Department of Statistics, University of California,
Berkeley, 367
Evans Hall, Berkeley, CA~94720-3860, USA. \printead{e3}}
\end{aug}

\received{\smonth{6} \syear{2011}}

%
\begin{abstract}
We develop a new formulation of Stein's method to obtain computable
upper bounds
on the total variation distance between the geometric distribution and a
distribution of interest. Our framework reduces the problem to the construction
of a coupling between the original distribution and the ``discrete equilibrium''
distribution from renewal theory. We illustrate the approach in four
non-trivial examples: the geometric sum of independent, non-negative,
integer-valued random variables having common mean, the generation size
of the
critical Galton--Watson process conditioned on non-extinction, the
in-degree of a
randomly chosen node in the uniform attachment random graph model and
the total
degree of both a fixed and randomly chosen node in the preferential attachment
random graph model.
\end{abstract}

%

\end{frontmatter}

\section{Introduction}\label{intr}

The exponential and geometric distributions are convenient and accurate
approximations in a wide variety of complex settings involving rare events,
extremes and waiting times.
The difficulty in obtaining explicit error bounds
for these approximations beyond
elementary settings is discussed in the
preface of~\cite{Aldous1989}, where the author also points out a
lack of such
results. Recently, Pek{\"o}z and R{\"o}llin~\cite{Pekoz2009} developed
a framework to obtain error
bounds
for the Kolmogorov and Wasserstein distance metrics between the exponential
distribution and a distribution of interest. The main ingredients there are
Stein's method (see~\cite{Ross2011,Ross2007} for introductions) along
with the
equilibrium distribution from renewal theory. Due to the flexibility of Stein's
method and the close connection between the exponential and geometric
distributions, it is natural to attempt to use similar techniques to obtain
bounds for the stronger total variation distance metric between the geometric
distribution and an integer supported distribution.

Our formulation rests on the idea that a positive, integer-valued
random variable
$W$ will be approximately geometrically distributed with parameter
$p=1/\IE W$
if
$\law(W)\approx\law(W^e)$, where~$W^e$ has the (discrete) equilibrium
distribution with respect to $W$, defined by
%
\begin{equation} \label{1}
\IP(W^e = k) = \frac{1}{\IE W} \IP(W\geq k),\qquad k=1,2,\ldots.
\end{equation}
This distribution arises in discrete-time renewal theory as the
time until the next renewal when the process is stationary, and the
transformation which maps a distribution to its equilibrium
distribution has
the geometric distribution with positive support
as its unique fixed point. Our main result is
an upper bound on the variation distance between the distribution of
$W$ and a
geometric distribution with parameter $(\IE W)^{-1}$, in terms of a coupling
between the
random variables $W^e$ and $W$.

This setup is closely related to the exponential approximation
formulation of
\cite{Pekoz2009} and also~\cite
{GoldsteinPC}, which is also related to
the zero-bias transformation of \cite
{Goldstein1997}. A difficulty in pushing
the results of~\cite{Pekoz2009} through to the
stronger total variation metric
is that the support of the
distribution to be approximated may not match the support of the geometric
distribution well enough.
This
issue is typical in bounding the total variation distance between integer-valued
random variables
and can be handled by introducing
a term into the bound that quantifies the ``smoothness''
of the distribution of interest; see, for
example,
\cite{Barbour2002,Rollin2005,Rollin2008a}. To illustrate this point, we apply
our abstract formulation to obtain total variation error bounds in
two of the examples treated in~\cite{Pekoz2009}.

It is also important to note that geometric approximation can be appropriate
in situations where exponential approximation is not; for example
if a sequence of random variables has the geometric distribution with fixed
parameter as its distributional limit.
Thus, we will
also apply our
theory in two examples (discussed in more detail immediately below)
that fall
into this category. In these applications, the ``smoothness''
term in our bounds does not play a part. The same effect can be
observed when
comparing translated Poisson and Poisson approximation; see~\cite{Rollin2005}
and~\cite{Barbour1992}.

The first application in this article is a bound on the total variation distance
between the geometric distribution and the sum of a geometrically distributed
number of independent, non-negative, integer-valued random variables
with common
mean. The distribution of such geometric convolutions have been
considered in
many places in the literature in the setting of exponential
approximation and
convergence; the book-length treatment is given in \cite
{Kalashnikov1997}. The
second application is a variation on the classical theorem of Yaglom
\cite{Yaglom1947},
describing the asymptotic behavior of the generation size of a critical
Galton--Watson process conditioned on non-extinction. This theorem has
a large
literature of extensions and embellishments; see, for example,
\cite{Lalley2009}.
Pek{\"o}z and R{\"o}llin~\cite{Pekoz2009} obtained a rate of
convergence for the Kolmogorov distance
between the generation size of a critical Galton--Watson process
conditioned on
non-extinction and the exponential distribution. Here we obtain an analogous
bound for the geometric distribution in total variation distance. The third
application is to the in-degree of a randomly chosen node in the uniform
attachment random graph discussed in \cite
{Bollobas2001}, and the final
application is to the total degree of both a fixed and a randomly
chosen node in
the preferential attachment random graph discussed in \cite
{Bollobas2001}. As
mentioned before, these examples do not derive from an exponential
approximation result.

Finally, we mention that there are other formulations of geometric approximation
using Stein's method. For example, Pek{\"o}z~\cite{Pekoz1996} and
Barbour and Gr{\"u}bel~\cite{Barbour1995a}
use the intuition that a positive, integer-valued random variable $W$
approximately has a geometric distribution with parameter $p=\IP(W=1)$ if
%
\[
\law(W) \approx\law(W-1|W>1).
\]
Other approaches can be found in
\cite{Phillips2000} and~\cite{Daly2010}.

The organization of this article is as follows. In Section~\ref{sec1}
we present
our main
theorems, and Sections~\ref{sec2},~\ref{sec3},~\ref{sec4},
and~\ref{sec5},
respectively, contain applications to geometric
sums, the critical Galton--Watson process conditioned on
non-extinction, the
uniform attachment random graph model and the
preferential attachment random graph model.

\section{Main results}\label{sec1}

A typical issue when discussing the geometric distribution is whether
to have the support begin at zero or one.
Denote by $\Ge(p)$ the geometric distribution with positive support;
that is,
$\law(Z)=\Ge(p)$ if $\IP(Z=k)=(1-p)^{k-1}p$ for positive integers $k$.
Alternatively, denote
by $\Ge^0(p)$ the geometric distribution $\Ge(p)$ shifted by minus one,
that is,
``starting at 0.'' Since $\law(Z)=\Ge(p)$ implies $\law(Z-1)=\Ge^0(p)$,
it is
typical that results for one of $\Ge(p)$ or $\Ge^0(p)$ easily pass to
the other.
Unfortunately, our methods do not appear to trivially transfer
between these two distributions, so we are forced to develop our theory
for both cases in parallel.

First, we give an alternate definition of the equilibrium distribution
that we
will use in the proof of our main result.

\begin{definition} \label{2} Let $X$ be a positive, integer-valued random
variable with finite mean. We say that an integer-valued random variable
$X^{e}$ has the
\textit{discrete equilibrium distribution w.r.t. $X$} if, for all
bounded $f$ and
$\nabla f(x) = f(x)-f(x-1)$, we have
%
\begin{equation} \label{3}
\IE f(X)-f(0) = \IE X \IE\nabla f(X^e).
\end{equation}
\end{definition}

\begin{remark}
To see how \eq{3} is equivalent to~\eq{1}, note that we have
\[
\IE f(X)-f(0)
= \IE\sum_{i=1}^X \nabla f(i)
= \sum_{i=1}^\infty\nabla f(i) \IP(X\geq i)
= \IE X \IE\nabla f(X^e).
\]
\end{remark}


Note that Definition~\ref{2} and \eq{1}
do not actually require positive support of
$W$ so that we can define $W^e$ for a non-negative random
variable which is not identically zero (and our results below still
hold in this case).
However, we will use $W^e$ to compare the distribution
of $W$ to a geometric distribution with positive support,
so the assumption that $W>0$ is not restrictive in most cases of interest.
In order to handle geometric
approximation with support on the non-negative integers, we could shift
the non-negative random variable
$W$ by 1 and then consider geometric approximation on the positive integers,
but this strategy is not practical
since $(W+1)^e$ is typically inconvenient to work with.
Fortunately, developing an analogous theory for
non-negative random variables with mass at zero is
no more difficult than that for positive random variables.

\begin{definition} If $X$ is a non-negative, integer-valued random variable
with $\IP(X=0)>0$, we say that an
integer-valued random variable $X^{e_0}$ has the \textit{discrete equilibrium
distribution w.r.t. $X$} if, for all bounded $f$ and with $\D f(x) =
f(x+1)-f(x)$, we have
%
\[
\IE f(X)-f(0) = \IE X \IE\D f(X^{e_0}).
\]
\end{definition}

It is not difficult the see that $W^{e_0} \eqlaw W^e -1$. Note that we are
defining the term ``discrete equilibrium distribution'' in both of the previous
definitions, but this should not cause confusion as the support of the base
distribution dictates the meaning of the terminology.

Besides the total variation metric, we will also give bounds on the
local metric
%
\[
\dloc(\law(U),\law(V)) := \sup_{m\in\IZ} |\IP(U=m)-\IP(V=m)|.
\]
It is clear that $\dloc$ will be less than or equal to
$\sup_{m}{[\IP(U=m)\vee\IP(V=m)]}$, so that typically better rates
need to be obtained in order to provide useful information
in this metric.

As a final bit of notation, before the statement of our main results,
for a function $g$ with domain~$\IZ$, let $\Vert{g}\Vert=\sup_{k\in\IZ}\vert{g(k)}\vert$,
and
for any
integer-valued random variable $W$ and any $\sigma$-algebra $\mathcal{F}$,
define the
conditional smoothness
%
\begin{equation} \label{4}
S_1(W|\mathcal{F})
= \sup_{\Vert{g}\Vert\leq1}\vert{\IE\{\Delta g(W)|\mathcal{F}\}}\vert
=2\dtv\bigl(\law(W+1|\mathcal{F}), \law(W|\mathcal{F})\bigr),
\end{equation}
and the second order conditional smoothness
\[
S_2(W|\mathcal{F})
 = \sup_{\Vert{g}\Vert\leq1}\vert{\IE\{\Delta^2 g(W)|\mathcal{F}\}}\vert,
\]
where $\Delta^2 g(k) = \Delta g(k+1)- \Delta g(k)$. In order to
simplify the
presentation of the main theorems, we let $d_1 = \dtv$ and $d_2 =
\dloc$.
Also, let $I_A$ denote the
indicator random variable of an event $A$.

\begin{theorem}\label{thm1} Let $W$ be a positive, integer-valued
random variable
with $\IE W=1/p$ for some $0<p\leq1$, and let $W^e$ have the discrete
equilibrium
distribution
w.r.t.~$W$. Then, with $D=W-W^e$, any $\sigma$-algebra $\mathcal{F}
\supseteq
\sigma(D) $ and event $A\in\mathcal{F}$, we have
\begin{equation}
d_l(\law(W),\Ge(p))
\leq\IE\{|D|S_l(W|\mathcal{F})I_A\}+2\IP(A^c)\label{5}
\end{equation}
for $l=1,2$, and
\begin{eqnarray}
\dtv(\law(W^e),\Ge(p))&\leq& p\IE\vert{D}\vert, \label{6}\\
\dloc(\law(W^e),\Ge(p))&\leq& p\IE\{\vert{D}\vert S_1(W|\mathcal{F})\};
\label{7}
\end{eqnarray}
on the RHS of~\eq{5} and~\eq{7}, $S_{l}(W|\mathcal{F})$ can be replaced by
$S_{l}(W^e|\mathcal{F})$.
\end{theorem}

\begin{theorem}\label{thm2} Let $W$ be a non-negative, integer-valued random
variable with $\IP(W=0)>0$, $\IE W=(1-p)/p$ for some $0<p\leq1$, and let
$W^{e_0}$ have the
discrete equilibrium distribution w.r.t.~$W$. Then, with $D=W-W^{e_0}$, any
$\sigma$-algebra $\mathcal{F} \supseteq\sigma(D) $ and event $A\in
\mathcal{F}$,
we have
\begin{equation}
d_l(\law(W),\Ge^0(p))
\leq(1-p)\IE\{|D|S_l(W|\mathcal{F})I_A\}+2(1-p)\IP(A^c) \label{8}
\end{equation}
for $l=1,2$, and
\begin{eqnarray}
\dtv(\law(W^{e_0}),\Ge^0(p))&\leq& p\IE\vert{D}\vert, \label{9}\\
\dloc(\law(W^{e_0}),\Ge^0(p))&\leq&
p\IE\{\vert{D}\vert S_1(W|\mathcal{F})\}, \label{10}
\end{eqnarray}
on the RHS of~\eq{8} and~\eq{10}, $S_l(W|\mathcal{F})$ can be replaced by
$S_l(W^{e_0}|\mathcal{F})$.
\end{theorem}

Before we prove Theorems~\ref{thm1} and~\ref{thm2}, we make a few remarks
related to these results.

\begin{remark}\label{rem2}
It is easy to see that the a random variable $W$ with law equal to $\Ge
(p)$ has
the property that $\law(W)=\law(W^e)$, so that $W^e$ can be taken to be
$W$ and
the theorem yields the correct error term in this case. The analogous statement
is true for $\Ge^0(p)$ and~$W^{e_0}$.
\end{remark}

\begin{remark}\label{rem3}
By choosing $A=\{W=W^{e_0}\}$ in \eq{8} with $l=1$, we find
\begin{equation}
\dtv(\law(W),\Ge^0(p)) \leq2(1-p)\IP(W\not=W^{e_0}), \label{11}
\end{equation}
and an analogous corollary holds for Theorem~\ref{thm1}.
\end{remark}

In order to use the theorem we need to be able to construct random variables
with the discrete equilibrium
distribution.
The next proposition provides such a construction for a non-negative,
integer-valued random variable $W$. We say
$W^s$ has the size-bias distribution of $W$, if
%
\[
\IE\{W f(W)\} = \IE W\IE f(W^s)
\]
for all $f$ for which the expectation exists.

\begin{proposition} \label{12}
Let $W$ be an integer-valued random variable, and let $W^s$ have
the size-bias distribution of $W$.
\begin{enumerate}[1.]
\item[1.] If $W>0$ and we define the random variable $W^e$
such that conditional on $W^s$, $W^e$ has the uniform distribution on the
integers $\{1, 2, \ldots, W^s\}$, then $W^e$ has the discrete equilibrium
distribution w.r.t.~$W$.
\item[2.] If $W\geq0$ with $\IP(W=0)>0$, and we define the random variable $W^{e_0}$
such that conditional on $W^s$, $W^{e_0}$ has the uniform distribution
on the
integers $\{0, 1, \ldots,\break W^s-1\}$, then $W^{e_0}$ has the discrete equilibrium
distribution w.r.t.~$W$.
\end{enumerate}
\end{proposition}
\begin{pf}
For any bounded $f$ we have
\[
\IE f(W)-f(0)
= \IE\sum_{i=1}^{W} \nabla f(i)
= \IE W\IE\Biggl\{\frac{1}{W^s}\sum_{i=1}^{W^s} \nabla f(i)\Biggr\}
= \IE W\IE\nabla f(W^e),
\]
which implies Item 1. The second item is proved analogously.
\end{pf}

As mentioned in the \hyperref[intr]{Introduction}, there can be considerable technical difficulty
in ensuring the support of the distribution to be approximated is smooth.
In Theorems~\ref{thm1} and~\ref{thm2} this issue is accounted for in
the term
$S_1(W|\mathcal{F})$. Typically, our strategy to bound this term will
be to
write $W$
(or $W^e$) as a sum of terms which are independent given $\mathcal{F}$ and then apply
the following lemma from~\cite{Mattner2007}, Corollary~1.6.
\begin{lemma}[(\cite{Mattner2007}, Corollary~1.6)]
\label{lem1}
If $X_1,\ldots, X_n$ are independent, integer-valued random variables and
%
\[
u_i = 1-\dtv\bigl(\law(X_i),\law(X_i+1)\bigr),
\]
then
%
\[
\dtv\Biggl(\law\Biggl(\sum_{i=1}^n X_i\Biggr),
\law\Biggl(1+\sum_{i=1}^n X_i\Biggr)\Biggr)
\leq
\sqrt{\frac{2}{\uppi}} \Biggl(\frac{1}{4} + \sum_{i=1}^n
u_i\Biggr)^{-1/2}.
\]
\end{lemma}


Before we present the proof of Theorems~\ref{thm1} and~\ref{thm2},
we must first develop the Stein method machinery we will need.
As in~\cite{Pekoz1996}, for any subset $B$ of the integers
and any
$p=1-q$, we construct the function $f = f_{B,p}$ defined by $f(0)=0$
and for
$k\geq1$,
\begin{equation}
qf(k)-f(k-1)=I_{k\in B} -\Ge(p)\{B\}, \label{13}
\end{equation}
where $\Ge(p)\{B\}=\sum_{i\in B} (1-p)^{i-1}p$ is the chance that
a positive, geometric random variable with parameter $p$
takes a value in the set $B$.
It can be easily verified that the
solution of \eq{13} is given by
%
\begin{equation} \label{14}
f(k)=\sum_{i\in B }q^{i-1} - \sum_{i\in B, i\geq k+1} q^{i-k-1}.
\end{equation}
Equivalently, for $k\geq0$,
%
\[
qf(k+1)-f(k)=I_{k\in B-1} -\Ge^0(p)\{B-1\},
\]
where we define $\Ge^0(p)\{B\}$ analogously to $\Ge(p)\{B\}$.\vadjust{\goodbreak}

Pek{\"o}z~\cite{Pekoz1996} and Daly~\cite{Daly2008} study properties
of these
solutions, but we need the following additional lemma to obtain our main
result.

\begin{lemma} \label{lem2} For $f$ as above,
we have
%
\begin{equation} \label{15}
\sup_{k\geq1}|\nabla f(k)| = \sup_{k\geq0}|\D f(k)|\leq1.
\end{equation}
If, in addition, $B=\{m\}$ for some $m\in\IZ$, then
%
\begin{equation} \label{16}
\sup_{k\geq0}\vert{f(k)}\vert \leq1.
\end{equation}
\end{lemma}
\begin{pf} To show \eq{15}, note that
\begin{eqnarray*}
\nabla f(k)
& = &\sum_{i\in B, i\geq k} q^{i-k}-\sum_{i\in B, i\geq k+1} q^{i-k-1}
\\
& =& I_{k\in B}+\sum_{i\in B, i\geq k+1} (q^{i-k}- q^{i-k-1})= I_{k\in B}-p \sum_{i\in B, i\geq k+1} q^{i-k-1},
\end{eqnarray*}
thus $-1\leq\nabla f(k) \leq1$. If now $B=\{m\}$, \eq{16} is
immediate from
\eq{14}.
\end{pf}

\begin{pf*}{Proof of Theorem~\protect\ref{thm1}}
Given any positive, integer-valued random variable $W$ with $\IE W=1/p$ and
$D=W-W^e$ we have, using \eq{13}, Definition~\ref{2}, and Lemma~\ref
{lem2} in
the two inequalities,
\begin{eqnarray*}
&&\IP(W\in B)-\Ge(p)\{B\}\\
& &\quad= \IE\{qf(W)-f(W-1)\}\\
&&\quad = \IE\{\nabla f(W)-pf(W)\}\\
&&\quad = \IE\{\nabla f(W)-\nabla f(W^e)\}\\
&&\quad \leq\IE\bigl\{I_A\bigl(\nabla f(W)-\nabla f(W^e)\bigr)\bigr\}+2\IP(A^c)\\
&&\quad = \IE\Biggl\{I_A I_{D>0}\sum_{i=0}^{D-1} \IE\bigl(\nabla
f(W^e+i+1)-\nabla f(W^e+i)|\mathcal{F}\bigr)\Biggr\}\\
& &\qquad{}+ \IE\Biggl\{I_A I_{D<0}\sum_{i=0}^{-D-1} \IE\bigl(\nabla
f(W^e-i-1)-\nabla f(W^e-i)|\mathcal{F}\bigr)\Biggr\}+2\IP(A^c)\\
&&\quad \leq\IE\{|D|S_1(W^e|\mathcal{F})I_A\}+2\IP(A^c),
\end{eqnarray*}
which is~\eq{5} for $l=1$; analogously, one can obtain~\eq{5} with
$S_1(W|\mathcal{F})$ in place of $S_1(W^{e}|\mathcal{F})$
on the RHS. In the case of $B=\{m\}$ we can make use of \eq{16} to obtain
%
\[
\vert{\IP(W=m)-\Ge(p)\{m\}}\vert \leq\IE\{|D|S_2(W^e|\mathcal{F})\I_A\}+2\IP(A^c),
\]
instead, which proves \eq{5} for $l=2$. For~\eq{6}, we have that
\begin{eqnarray*}
\IP(W^e\in B)-\Ge(p)\{B\}
&= &\IE\{qf(W^e)-f(W^e-1)\} \\
& =& \IE\{\nabla f(W^e)-pf(W^e)\}
= p\IE\{f(W)-f(W^e)\} \leq p\IE|D|,
\end{eqnarray*}
where the last line follows by writing $f(W)-f(W^e)$ as a telescoping
sum of
$|D|$ terms no greater than $\Vert{\nabla f}\Vert$, which can be bounded
using~\eq{15}; \eq{7} is straightforward using \eq{16} and~\eq{4}.
\end{pf*}

\begin{pf*}{Proof of Theorem~\protect\ref{thm2}}
Let $W$ be a non-negative, integer-valued random variable with $\IE W=(1-p)/p$
and $W^{e_0}$ as in the theorem.
If we define $I$ to be independent of all else and such that
$\IP(I=1)=1-\IP(I=0)=p$, then a short calculation shows that
the variable defined by
%
\[
[(W+1)^e| I=1]=W+1\quad \mbox{and} \quad[(W+1)^e|
I=0]=W^{e_0}+1
\]
has the positive discrete equilibrium transform with respect to $W+1$. Equation
\eq{8} now follows after noting that
%
\[
d_l(\law(W),\Ge^0(p))=
d_l\bigl(\law(W+1),\Ge(p)\bigr),
\]
and then applying the following stronger variation of \eq{5} which is easily
read from the proof of Theorem~\ref{thm1}:
\begin{eqnarray*}
 &&d_l\bigl(\law(W+1),\Ge(p)\bigr)\\
 &&\quad\leq\IE\{|W+1-(W+1)^e|S_l(W|\mathcal{F})I_A\}
+2\IP\bigl((W+1)^e\not=W+1\bigr)\IP(A^c).
\end{eqnarray*}
Equations \eq{9} and~\eq{10} can be proved in a manner similar to their analogs
in Theorem~\ref{thm1}.
\end{pf*}

\section{Applications to geometric sums}\label{sec2}

In this section we apply the results above to a sum of the geometric
number of
independent but not necessarily identically distributed random
variables. As in
our theory above, we will have separate results for the two cases where
the sum is strictly positive and the case where it can take on the
value zero
with positive probability.
We reiterate that although there are a variety of exponential approximation
results in the literature for this example,
there do not appear to be bounds available for the analogous geometric
approximation in the total variation metric.
\begin{theorem}\label{thm3} Let
$X_1, X_2, \ldots$ be a sequence of independent, square integrable,
positive
and integer-valued random variables, such
that, for some $u>0$, we have, for all $i\geq1$,
$\IE X_i = \mu$ and $u \leq1-\dtv(\law(X_i), \law(X_i+1))$.\vadjust{\goodbreak}
Let $\law(N)=\Ge(a)$ for some $0<a\leq1$ and $W =
\sum_{i=1}^N X_i$. Then with $p=1-q=a/\mu$, we have
\begin{equation}\label{17}
d_l(\law(W),\Ge(p))
\leq C_l \sup_{i\geq1} \mathbb{E} |X_i-X_i^e| \leq C_l
\biggl(\mu_2/2+\ahalf+\mu\biggr)
\end{equation}
for $l=1,2$, where $\mu_2: =\sup_i \IE X_i^2$ and
\begin{eqnarray*}
C_1 & =&\min\biggl\{1, a\biggl[1+\biggl(-\frac{2}{u\log
(1-a)}\biggr)^{1/2}\biggr]\biggr\},\\
C_2 & = &\min\biggl\{1, a\biggl[1-\frac{6\log(a)}{\uppi u}\biggr]\biggr\}.
\end{eqnarray*}
\end{theorem}
\begin{theorem}\label{thm4} Let
$X_1, X_2, \ldots$ be a sequence of independent, square integrable,
non-negative and integer-valued random variables, such
that for some $u>0$ we have, for all $i\geq1$,
$\IE X_i = \mu$ and $u \leq1-\dtv(\law(X_i), \law(X_i+1))$.
Let $\law(M)=\Ge^0(a)$ for some $0<a\leq1$ and $W =
\sum_{i=1}^M X_i$. Then with $p=1-q=a/(a+\mu(1-a))$, we have
%
\begin{equation}\label{18}
d_l(\law(W),\Ge^0(p))
\leq C_l \sup_{i\geq1} \mathbb{E} X_i^{e_0} \leq
C_l\biggl(\mu_2/(2\mu)-\ahalf\biggr)
\end{equation}
for $l=1,2$, where $\mu_2: =\sup_i \IE X_i^2$ and the $C_l$ are as in
Theorem~\ref{thm3}.
\end{theorem}
\begin{remark} The first inequality in \eq{17} yields the correct bound of zero when $X_i$
is geometric, since in this case we would have $X_i=X_i^e$; see
Remark~\ref{rem2}
following Theorem~\ref{thm1}. Similarly, in the case where the $X_i$
have a
Bernoulli
distribution with expectation $\mu$, we have that $X^{e_0} = 0$ so that
the left-hand side of \eq{18} is zero.
That is, if $M\sim\Ge^{0}(a)$ and conditional on $M$, $W$ has the binomial
distribution with
parameters $M$ and $\mu$ for some $0\leq\mu\leq1$, then
$W\sim\Ge^{0}(a/(a+\mu(1-a)))$.
\end{remark}
%
%
\begin{remark} In the case where $X_i$ are i.i.d. but not
necessarily integer
valued and $0<a\leq\ahalf$, Brown~\cite{Brown1990}, Theorem 2.1,
obtains the
exponential approximation result
\begin{equation}\label{19}
\dk\bigl(\law(W),\Exp(1/p)\bigr)
\leq\frac{a\mu_2}{\mu}
\end{equation}
for the weaker Kolmogorov metric.
To compare \eq{19} with \eq{17} for small $a$, we observe that bound \eq{17}
is linear in $a$ whereas \eq{19}, within a constant factor, behaves
like $a
(-\log(1-a))^{-1/2}\sim\sqrt{a}$. Therefore bound \eq{19} is better, but
\eq{17} applies to non-i.i.d. random variables (albeit having identical means)
and to the stronger total variation metric.
\end{remark}

\begin{pf*}{Proof of Theorem~\protect\ref{thm3}}
First, let us prove that $W^e := \sum_{i=1}^{N-1} X_i + X_N^e$ has the discrete
equilibrium distribution w.r.t. $W,$ where, for each $i\geq1$, $X_i^e$
is a random variable having the equilibrium distribution w.r.t. $X_i$,\vadjust{\goodbreak}
independent of all else. Note first
that we have, for bounded $f$ and every~$m$,
%
\[
\mu\IE\nabla f\Biggl(\sum_{i=1}^{m-1} X_i +X_m^e\Biggr)
= \IE\Biggl[ f\Biggl(\sum_{i=1}^m X_i\Biggr)
- f\Biggl(\sum_{i=1}^{m-1} X_i\Biggr)\Biggr].
\]
Note also that since $N$ is geometric, for any bounded function $g$ with
$g(0)=0$, we have $\IE\{{g(N)} - {g(N-1)}\} = a\IE g(N)$.
We now assume that $f(0) = 0$. Hence, using the above two facts and
independence between $N$ and the sequence $X_1,X_2,\ldots,$ we have
\begin{eqnarray*}
\IE W \IE\nabla f(W^e)
& =& \frac{\mu}{a} \IE\nabla f\Biggl(\sum_{i=1}^{N-1} X_i + X_N^e\Biggr)\\
& =& \frac{1}{a}\IE f\Biggl[\Biggl(\sum_{i=1}^{N} X_i\Biggr)-
f\Biggl(\sum_{i=1}^{N-1} X_i\Biggr)\Biggr] = \IE f\Biggl(\sum_{i=1}^{N} X_i\Biggr) = \IE f(W).
\end{eqnarray*}
Now, $D= W - W^e = X_N-X_N^e$ and setting $\mathcal{F}=\sigma(N, X_N^e, X_N)$,
we have
\begin{eqnarray} \label{20}
S_1(W^e|\mathcal{F}) &=& S_1\Biggl(\sum_{i=1}^{N-1} X_i \Big\bbbmid\mathcal{F}\Biggr)
 \leq1\wedge\biggl(\frac{2}{\uppi(0.25 + (N-1)u)}\biggr)^{1/2}
 \nonumber
 \\[-8pt]
 \\[-8pt]
 \nonumber
& \leq&1\wedge\biggl(\frac{2}{\uppi(N-1)u}\biggr)^{1/2},
\end{eqnarray}
where we have used Lemma~\ref{lem1} and the fact that $S_1(W^e|\mathcal
{F})$ is
almost surely bounded by one.
We now have
%
\begin{equation}
\IE[|D|S_1(W^e|\mathcal{F})] \leq\IE\biggl[
{\biggl(1\wedge{\biggl(\frac{2}{\uppi(N-1) u}\biggr)}^{1/2}\biggr)} \IE^N | X_N-X^e_N
|\biggr]. \label{21}
\end{equation}

From here, we can obtain the first inequality in~\eq{17} by applying
Theorem~\ref{thm1}, after noting that
\begin{equation}
\IE^N | X_N-X^e_N|\leq\sup_{i\geq1}\IE\vert{X_i-X^e_i}\vert,
\end{equation}
and
\begin{eqnarray*}
\IE{\biggl(1\wedge{\biggl(\frac{2}{\uppi(N-1) u}\biggr)}^{1/2}\biggr)}
& \leq& 1\wedge\biggl(a + {\biggl(\frac{2}{\uppi u}\biggr)}^{1/2}\sum_{i\geq
1}\frac{a(1-a)^i}{i^{1/2}}\biggr)
\\
& \leq& 1\wedge\biggl(a + a{\biggl(\frac{2}{\uppi
u}\biggr)}^{1/2}{\biggl(-\frac{\uppi}{\log(1-a)}\biggr)}^{1/2}\biggr)\\
& =& 1\wedge\biggl(a{\biggl[1+ {\biggl(-\frac{2}{\log(1-a)u}\biggr)}^{1/2}\biggr]}\biggr),
\end{eqnarray*}
where we have used
%
\[
\sum_{i\geq1} \frac{a(1-a)^{i-1}}{i^{1/2}} \leq\frac{a}{1-a}
\int_0^\infty\frac{(1-a)^x}{x^{1/2}}\, \mathrm{d}x = \frac{a}{1-a}{\biggl(-\frac
{\uppi}{\log
(1-a)}\biggr)}^{1/2}.
\]
The second inequality in~\eq{17} follows from Theorem~\ref{thm1} and
the fact
(from the definition of the transformation $X^e$) that $\IE^N | X_N-X^e_N
| \leq\mu_2/2+\ahalf+\mu$.

To obtain the local limit result, note that,
if $V=X+Y$ is the sum of two independent random variables, then
$S_2(V)\leq
S_1(X)S_1(Y)$. Hence,
\[
S_2(W^e|\mathcal{F})
 \leq1\wedge\frac{2}{\uppi(0.25 + (N/2-1)_+u)}
\leq1\wedge\frac{6}{\uppi(N-1)u}.
\]
From here we have
%
\[
\IE S_2(W^e|\mathcal{F}) \leq1\wedge{\biggl(a+\frac{6}{\pi u}\sum_{i\geq
1}\frac{a(1-a)^i}{i}\biggr)} = 1\wedge{\biggl(a-\frac{6a\log(a)}{\uppi u}\biggr)}.
\]
\upqed\end{pf*}

\begin{pf*}{Proof of Theorem~\protect\ref{thm4}} It is straightforward to
check that
$W^{e_0} := \sum_{i=1}^M X_i + X_{M+1}^{e_0}$ has the equilibrium distribution
with respect to $W$.
Now, $D= W -W^{e_0} = -X_{M+1}^{e_0}$ and setting $\mathcal{F}=\sigma(M,
X_{M+1}^{e_0})$, we have
\[
S_1(W^{e_0}|\mathcal{F}) = S_1{\Biggl(\sum_{i=1}^{M} X_i \Big\bbbmid\mathcal{F}\Biggr)},
\]
which can be bounded above by \eq{20}, as in the proof of Theorem~\ref{thm3}.
The remainder of the proof
follows closely to that of Theorem~\ref{thm3}. For example, the expression
analogous to \eq{21} is
%
\[
\IE[|D|S_1(W^{e_0}|\mathcal{F})] \leq\IE\biggl[ \min\biggl\{1,
\frac{\sqrt{2}}{(\uppi M u)^{1/2}}\biggr \} \IE^M X^{e_0}_M
\biggr],
\]
and the definition of the transform $X^{e_0}$ implies that
%
\[
\IE X_i^{e_0}=\frac{\IE X_i^2}{2\mu}-\frac{1}{2}.
\]
\upqed\end{pf*}

\section{Application to the critical Galton--Watson branching
process}\label{sec3}

Let $Z_0=1,Z_1,Z_2,\ldots$ be a Galton--Watson branching process with offspring
distribution $\law(Z_1)$. A theorem due to Yaglom~\cite{Yaglom1947}
states that, if
$\IE Z_1 = 1$ and $\Var Z_1 = \sigma^2<\infty$, then $\law( n^{-1}Z_n |Z_n>0)$
converges to an exponential distribution with mean $\sigma^2/2$. The recent
article~\cite{Pekoz2009} is the first to give
an explicit bound on the rate of
convergence for this asymptotic result.\vadjust{\goodbreak} Using ideas from there, we give a
convergence rate for the total variation error of a geometric
approximation to
$Z_n$ under finite third moment of the offspring distribution and the
natural periodicity requirement that
%
\begin{equation} \label{22}
\dtv\bigl(\law(Z_1), \law(Z_1+1)\bigr)<1.
\end{equation}
This type of smoothness condition is typical in the context of Stein's
method for approximation by a discrete distribution; see, for example,
\cite{Barbour2002}
and~\cite{Rollin2008a}.

For the proof of the following theorem, we make
use the of construction of Lyons \textit{et~al.}~\cite{Lyons1995}; we
refer to
that article for more details on the construction and only present what is
needed for our purpose.

\begin{theorem}\label{thm5} For a critical Galton--Watson branching
process with
offspring distribution $\law(Z_1)$, such that $\IE Z_1^3 < \infty$
and \eq{22} hold, we have
%
\begin{equation}\label{23}
\dtv\biggl(\law(Z_n|Z_n>0),\Ge\biggl(\tsfrac{2}{\sigma^2n}\biggr)\biggr)
\leq\frac{C\log n}{n^{1/4}}
\end{equation}
for some constant $C$ which is independent of $n$.
\end{theorem}
\begin{remark}From~\cite{Pekoz2009}, Theorem 4.1, we have
%
\begin{equation}\label{24}
\dk\bigl(\law\bigl(2Z_n/(\sigma^2n)|Z_n>0\bigr),\Exp(1)\bigr)
\leq C\biggl(\frac{\log n}{n}\biggr)^{1/2}
\end{equation}
without condition \eq{22} for the weaker Kolmogorov metric. It can
be seen
that the bound in \eq{23} is not as good as the bound in \eq{24} for
large $n$,
but \eq{23} applies to the stronger total variation metric.
\end{remark}

\begin{pf*}{Proof of Theorem~\protect\ref{thm5}}
First we construct a size-biased branching tree as in~\cite{Lyons1995}. We
assume that this tree is labeled and ordered, in the sense that, if $w$
and $v$
are vertices in the tree from the same generation, and $w$ is to the
left of $v$,
then the offspring of $w$ is to the left of the offspring of $v$. Start in
generation $0$ with one vertex $v_0$, and let it have a number of offspring
distributed according to the size-bias distribution of $\law(Z_1)$.
Pick one of
the
offspring of $v_0$ uniformly at random, and call it $v_1$. To each of the
siblings of $v_1$ attach an independent Galton--Watson branching
process with
offspring distribution $\law(Z_1)$. For $v_1$ proceed as for $v_0$,
that is, give
it a
size-biased number of offspring, pick one uniformly at random, call it $v_2$,
attach independent Galton--Watson branching process to the siblings of
$v_2$ and
so on. It is clear that this will always give an infinite tree as the ``spine''
$v_0,v_1,v_2,\ldots$ is an infinite sequence where $v_i$ is an
individual (or
particle) in generation~$i$.

We now need some notation. Denote by $S_n$ the total number of
particles in
generation $n$. Denote by $L_n$ and $R_n$ the number of particles
to the left (excluding $v_n$) and to the right (including $v_n$),
of vertex $v_n$, so that $S_n = L_n + R_n$. We can describe these
particles in
more detail, according to the generation at which they split off from
the spine.
Denote by $S_{n,j}$ the number of particles in generation $n$ that stem
from any
of the siblings of $v_j$ (but not $v_j$ itself). Clearly,\vadjust{\goodbreak} $S_n = 1 +
\sum_{j=1}^n S_{n,j}$, where the summands are independent. Likewise, let
$L_{n,j}$ and $R_{n,j}$, be the number of particles in generation
$n$ that stem from the siblings to the left and right of $v_j$
(note that $L_{n,n}$ and $R_{n,n}$ are just the number of siblings of
$v_n$ to
the left and to the right, respectively). We have the relations $L_n =
\sum_{j=1}^n L_{n,j}$ and $R_n = 1 + \sum_{j=1}^n R_{n,j}$. Note that, for
fixed~$j$, $L_{n,j}$ and $R_{n,j}$ are, in general, not independent, as
they are
linked through the offspring size of $v_{j-1}$.

Let now $R_{n,j}'$ be independent random variables such that
%
\[
\law(R'_{n,j}) = \law(R_{n,j} | L_{n,j} = 0),\vspace*{-1pt}
\]
and, with $A_{n,j} = \{ L_{n,j}=0\}$, define
%
\[
R_{n,j}^* = R_{n,j} I_{A_{n,j}} + R_{n,j}' I_{A_{n,j}^c}
= R_{n,j} + (R_{n,j}' - R_{n,j}) I_{A_{n,j}^c}.
\]
Define also $R_n^* = 1 + \sum_{j=1}^n R_{n,j}^*$. Let us collect a few
facts from~\cite{Pekoz2009} which we will then
use to give the
proof of the theorem (here and in the rest of the proof, $C$ shall
denote a
constant which is independent of $n$, but may depend on $\law(Z_1)$ and
may also
be different from formula to formula):
\begin{enumerate}[(iii)]
\item[(i)]$\law(R_n^*) = \law(Z_n | Z_n > 0)$;

\item[(ii)] $S_n$ has the size-biased distribution of $Z_n$,
and $v_n$ is equally likely to be any of the $S_n$
particles;

\item[(iii)] $\IE\{R_{n,j}'I_{A_{n,j}^c}\}\leq\sigma^2\IP
[A_{n,j}^c];$

\item[(iv)] $\IE\{R_{n,j} I_{A_{n,j}^c}\} \leq\gamma
\IP[A_{n,j}^c],$ and $ \IE\{R_{n-1,j} I_{A_{n,j}^c}\} \leq\gamma
\IP[A_{n,j}^c]$,
where $\gamma= \IE Z_1^3$;

\item[(v)] $\IP[A_{n,j}^c]\leq\sigma^2\IP[Z_{n-j}>0]\leq C /
(n-j+1)$ for some $C>0$.
\end{enumerate}
In light of (i) and (ii) (and then using the construction in
Proposition~\ref{12}
to see that $R_n$ has the discrete equilibrium distribution w.r.t $\law(R_n^*)$)
we can let $W=R_n^*$, $W^e = R_n$ and\looseness=-1
\[
D=R_n^*-R_n= \sum_{j=1}^n(R_{n,j}' - R_{n,j}) I_{A_{n,j}^c}.\vspace*{-1pt}
\]\looseness=0
Also let
%
\[
N=\sum_{j=1}^{n-1} R_{n-1,j}I_{A^c_{n,j}}\quad \mbox{and}\quad M=\sum_{j=1}^{n}
R_{n,j}I_{A^c_{n,j}},\vspace*{-1pt}
\]
and note that (iii)--(v) give
%
\begin{equation}\label{25}
\IE|D| \leq C\log n\quad \mbox{and}\quad \IE N \leq C \log n.
\end{equation}

Next with
$\mathcal{F}=\sigma(N, D, R_{n-1}, M, R_{n,n}I_{A_{n,n}})$
and, letting $Z_1^i$, $i=1,2,\ldots,$ be i.i.d. copies of~$Z_1$, we have
%
\[
\law(R_{n}-M- R_{n,n}I_{A_{n,n}}-1| \mathcal{F}) =
\law{\Biggl( \sum_{i=1}^{R_{n-1}-N} Z_1^i \Big\bbbmid R_{n-1},N\Biggr)},
\]
which follows since $R_{n}-M=1+\sum_{i=1}^n R_{n,j}I_{A_{n,j}}$,
and the particles counted by $R_{n-1}-N$ will be parents of the particles
counted by $R_{n}-M-1+R_{n,n}I_{A_{n,n}}$.\vadjust{\goodbreak}

Then we use Lemma~\ref{lem1} to obtain
\begin{equation}\label{26}
S_1(W^e|\mathcal{F}) = S_1(R_{n}-M- R_{n,n}I_{A_{n,n}}-1| \mathcal{F})
\leq\frac{0.8}{(0.25 + (R_{n-1}-N)u)^{1/2}}.
\end{equation}
As a direct corollary of~\eq{6}, for any bounded function $f$, we have
%
\begin{equation}\label{27}
\IE f(W^e) \leq\IE f(X_p)+ p\Vert{f}\Vert\IE\vert{W^e-W}\vert,
\end{equation}
where $X_p\sim\Ge(p)$.
Fix $q= 1/\IE[Z_{n-1}|Z_{n-1}>0]$, $k=q^{-1/4}$, and let
\[
A=\{N\leq k,  |D|\leq k,  R_{n-1} > 2k\},
\]
and
%
\[
f(x) = (x-k)^{-1/2}I_{x\geq2k+1}.
\]
Using \eq{25}, \eq{26}, \eq{27} and the fact that $\Vert{f}\Vert\leq k^{-1/2}$,
we find
\[
\IE[f(X_q)] \leq q\sum_{j=1}^\infty\frac{(1-q)^j}{j^{1/2}} \leq(q\uppi)^{1/2}
\]
to obtain
\begin{eqnarray*}
\IE[|D|S_1(W^e|\mathcal{F})I_A]
& \leq& k u^{-1/2} \IE f(R_{n-1}) \\
& \leq& k u^{-1/2}\bigl( \IE f(X_q)+qk^{-1/2}\IE|D_{n-1}|\bigr) \\
& \leq& C q^{1/4}\log n,
\end{eqnarray*}
where $D_{n-1}=R_{n-1}-R_{n-1}^*$.
Now, applying~\eq{6} yields
\[
\IP(R_{n-1}\leq2k) \leq1-(1-q)^{2k}+q\IE|D_{n-1}|\leq q(2k+\IE|D_{n-1}|),
\]
and, by Markov's inequality and \eq{25}, we finally obtain
%
\[
\IP(A^c)\leq k^{-1}(\IE N +\IE|D|) +q(2k+\IE|D|)
\leq C q^{1/4}\log n.
\]
The theorem follows after using $(v)$ and $\IE Z_n =1$ to obtain $\IE
[Z_n|Z_n>0]
\leq C n$.
\end{pf*}

\section{Application to the uniform attachment random graph model}\label{sec4}

Let $G_n$ be a directed random graph on $n$ nodes defined by the following
recursive construction.
Initially the graph starts with one node with a single loop, where one
end of the loop contributes to the ``in-degree'' and the other to the
``out-degree.''
Now, for $2\leq m\leq n$, given the graph with $m-1$ nodes, add node
$m$ along with an edge directed from $m$ to a node chosen uniformly at random
among
the $m$ nodes present. Note that this model allows edges connecting a
node with
itself. This random graph model is referred to as uniform attachment.

This model has been well studied, and it was shown in \cite
{Bollobas2001} that
if $W$ is equal to the in-degree of a node chosen uniformly at random from
$G_n$, then $W$ converges to a geometric distribution (starting at 0) with
parameter $1/2$ as $n\rightarrow\infty.$ We will give an explicit
bound on
the total variation distance between the distribution of $W$ and the geometric
distribution that yields this asymptotic.

The same result, but with a larger constant, was obtained in \cite
{Ford2009},
where the author uses Stein's method an ad hoc analysis of the model.

\begin{theorem}\label{thm6} If $W$ is the in-degree of a node chosen uniformly
at random from the random graph $G_n$ generated according to uniform attachment,
then
\[
\dtv{\biggl(\law(W),\Ge^0\biggl(\ahalf\biggr)\biggr)}\leq\frac{1}{n}.
\]
\end{theorem}

\begin{pf}
Let $X_i$ have a Bernoulli
distribution, independent of all else, with parameter $\mu_i := (n-i+1)^{-1}$,
and let $N$ be an independent random variable that is uniform on the integers
$1,2,\ldots, n$. If we imagine that node $n+1-N$ is the randomly
selected node,
then it's easy to see that we can write $W:=\sum_{i=1}^{N} X_i.$

Next, let us prove that $ \sum_{i=1}^{N-1} X_i$ has the discrete
equilibrium distribution w.r.t. $W.$ First note
that we have, for bounded $f$ and every $m$,
%
\[
\mu_m \IE\D f{\Biggl( \sum_{i=1}^{m-1} X_i\Biggr)}
= \IE{\Biggl[f{\Biggl( \sum_{i=1}^m X_i\Biggr)}
- f\Biggl( \sum_{i=1}^{m-1} X_i\Biggr)\Biggr]},
\]
where we use
%
\[
\IE f(X_m)-f(0) = \IE X_m \IE\D f(0)
\]
and thus the fact that we can write $(X_m)^{e_0} \equiv0$. Note also
that, for
any bounded function $g$ with
$g(0)=0$, we have
%
\[
\IE{\Biggl(\frac{g(N)}{\mu_N} - \frac{g(N-1)}{\mu_N}\Biggr)}
= \IE g(N).
\]
We now assume that $f(0) = 0$. Hence, using the above two facts and
independence between $N$ and the sequence $X_1,X_2,\ldots,$ we have
\[
\IE W \IE\D f{\Biggl( \sum_{i=1}^{N-1} X_i \Biggr)}=\IE f(W).
\]
Now, let
%
\[
N' =
\cases{
N & \quad$\mbox{if $1\leq N < n$,}$\vspace*{2pt}\cr
0 & \quad$\mbox{if $N=n$.}$}
\]
We have that $\law(N')=\law(N-1)$ so that $W^{e_0}:=\sum_{i=1}^{N'} X_i$
has the equilibrium distribution with respect to $W$, and it is plain that
%
\[
\IP[W \neq W^{e_0} ] \leq\IP[N=n] = \frac{1}{n}.
\]
Applying \eq{11} of Remark~\ref{rem3} yields the theorem.\vspace*{-2pt}
\end{pf}

\section{Application to the preferential attachment random graph
model}\label{sec5}

Define the directed graph $G_n$ on $n$ nodes by the following recursive
construction. Initially the graph starts with one node with a single
loop where
one end of the loop contributes to the ``in-degree'' and the other to the
``out-degree.'' Now, for $2\leq m\leq n$, given the graph with $m-1$
nodes, add
node $m$ along with an edge directed from $m$ to a random node chosen
proportionally to the total degree of the node. Note that at step $m$,
the chance
that node $m$ connects to itself is $1/(2m-1)$, since we consider the added
vertex $m$ as immediately having out-degree equal to one. This random graph
model is referred to as preferential attachment.

This model has been well studied, and it was shown in \cite
{Bollobas2001} that
if $W$ is equal to the in-degree of a node chosen uniformly at random from
$G_n$, then $W$ converges to the Yule--Simon distribution (defined
below). We
will give a rate of convergence in the total variation distance for this
asymptotic, a result that cannot be read from the main results in
\cite{Bollobas2001}. Some rates of
convergence in this and related random graph
models can be found in the thesis~\cite{Ford2009}, but the
techniques and
results there do not appear to overlap with ours below. Detailed
asymptotics for
the individual degrees, along with rates, can be found
in~\cite{Pekoz2011a}.\looseness=-1

We say the random variable $Z$ has the Yule--Simon distribution if
%
\[
\IP(Z=k) = \frac{4}{k(k+1)(k+2)},\qquad k=1,2, \ldots.
\]
The following is our main result.\vspace*{-2pt}
\begin{theorem} \label{thm7}
Let $W_{n,i}$ be the total degree of vertex $i$ in the
preferential attachment graph on $n$ vertices, and let $I$
uniform on $\{1, \ldots, n\}$ independent of $W_{n,i}$. If $Z$ has the
Yule--Simon distribution, then
%
\[
\dtv{(\law(W_{n,I}), \law(Z))}\leq\frac{C \log n}{n}
\]
for some constant $C$ independent of $n$.
\end{theorem}
\begin{remark} The notation $\law(W_{n,I})$ in the statement of Theorem~\ref{thm7} should
be interpreted as
%
\[
\law(W_{n,I}|I=i)=\law(W_{n,i}).
\]
We will use similar notation in what follows without further
mention.\vadjust{\goodbreak}
\end{remark}

\begin{proposition}\label{prop1}
If $U$ has the uniform distribution on $(0,1)$, and given~$U$, we
define $Z$
such that $\law(Z)=\Ge(\sqrt{U})$, then $Z$ has the Yule--Simon
distribution.\vspace*{-2pt}
\end{proposition}

Our strategy to prove Theorem~\ref{thm7} will be to show that, for $I$
uniform on
$\{1, \ldots, n\}$ and~$U$ uniform on $(0,1)$, we have:
\begin{enumerate}[(1)]
\item[(1)]\label{(1)} $\dtv{(\law(W_{n,I}), \Ge(\IE
[W_{n,I}|I]^{-1}))}\leq
C\log(n)/n$,
\item[(2)]\label{(2)} $\dtv{(\Ge(\IE[W_{n,I}|I]^{-1}), \Ge(\sqrt
{I/n}))}\leq
C\log(n)/n$,
\item[(3)]\label{(3)} $\dtv{( \Ge(\sqrt{I/n}), \Ge(\sqrt{U}))}\leq
C\log(n)/n$,
\end{enumerate}
where, here and in what
follows, we use the letter $C$ as a generic constant which may differ
from line
to line.
From this point, Theorem~\ref{thm7} follows from the triangle
inequality and
Proposition~\ref{prop1}.

Item~(1) will follow from our framework above; in particular we
will use the following result, which may be of independent interest. We postpone
the proof to the end of the section.\vspace*{-2pt}
\begin{theorem} \label{thm8}
Retaining the notation and definitions above, we have
%
\[
\dtv{\bigl(\law(W_{n,i}), \Ge\bigl(1/\IE(W_{n,i})\bigr)\bigr)}
\leq\frac{C}{i}
\]
for some constant $C$ independent of $n$ and $i$.\vspace*{-2pt}
\end{theorem}

To show Items~(2) and~(3) we will need the
following lemma. The first statement is found in~\cite{Bollobas2001}, page~283,
and the second
follows easily from the first.\vspace*{-2pt}
\begin{lemma}[(Bollob{\'a}s \textit{et~al.}~\cite{Bollobas2001})]\label{lem3}
Retaining the notation and definitions above, for all $1\leq i\leq n$,
%
\[
\Biggl|\IE W_{n,i} -\sqrt{\frac{n}{i}}\Biggr| \leq
C\sqrt{\frac{n}{i^{3}}}
\quad\mbox{and}\quad
\Biggl|\frac{1}{\IE W_{n,i}} -\sqrt{\frac{i}{n}}\Biggr| \leq
\frac{C}{\sqrt{ni}}.\vspace*{-2pt}
\]
\end{lemma}

Our final general lemma is useful for handling total variation distance
for conditionally defined random variables.\vspace*{-2pt}

\begin{lemma}\label{lem4} Let $W$ and $V$ be random variables, and let
$X$ be a
random element defined on the same probability space. Then
%
\[
\dtv(\law(W), \law(V))\leq\IE\dtv(\law(W|X), \law(V|X)).\vspace*{-2pt}
\]
\end{lemma}
\begin{pf} If $f\dvtx \IR\rightarrow[0,1]$, then
%
\[
\vert{\IE[f(W)-f(V)]}\vert\leq\IE\vert{\IE[f(W)-f(V)|X]}\vert
\leq\IE\dtv{(\law(W|X), \law(V|X))}.\vspace*{-2pt}
\]
\upqed\end{pf}

\begin{pf*}{Proof of Theorem~\protect\ref{thm7}}
We first claim that
%
\begin{equation}
\dtv{\bigl(\Ge(p),\Ge(p-\eps)\bigr)}
\leq\frac{\eps}{p}< \frac{\eps}{p-\eps}. \label{28}\vadjust{\goodbreak}
\end{equation}
The second inequality of \eq{28} is obvious. To see the first inequality,
we construct two infinite sequences of independent random variables.
The first sequence consists of $\Be(p)$ random variables, and the second
sequence consists of $\Be(p-\eps)$ random variables maximally
coupled component-wise to the first so that the
terms in the first sequence are no smaller than the
corresponding terms in the second. For each of these
sequences, the index of the first Bernoulli random variable, which is
$1$, follows
a $\Ge(p)$ and $\Ge(p-\eps)$ distribution,
respectively. Since the index of the first occurrence of a
$1$ in the
first sequence is less than or equal
to that in the second sequence,
the probability that these two random variables are not equal is the
probability that a
coordinate in the second sequence is $0$, given the same coordinate
is $1$ in the first
sequence, which is $\eps/p$.

Using \eq{28} and Lemma~\ref{lem3} we easily obtain
%
\[
\dtv{\bigl(\Ge(1/\IE W_{n,i}), \Ge\bigl(\sqrt{i/n}\bigr)\bigr)}
\leq\frac{C}{i},
\]
and applying Lemma~\ref{lem4}
we find
%
\[
\dtv{\bigl(\Ge(\IE[W_{n,I}|I]^{-1}), \Ge\bigl(\sqrt{I/n}\bigr)\bigr)}
\leq\frac{C\log(n)}n,
\]
which is Item~(2) above. Now, coupling $U$ to $I$ by writing $U=I/n-V$,
where $V$ is uniform on $(0,1/n)$
and independent of $I$,
and using first \eq{28} and then Lemma~\ref{lem4} leads to
%
\[
\dtv{\bigl(\Ge\bigl(\sqrt{U}\bigr), \Ge\bigl(\sqrt{I/n}\bigr)\bigr)}
\leq\frac{C}{n}\sum_{i=1}^n \frac{
\sqrt{i/n}-\sqrt{(i-1)/n}}{\sqrt{i/n}}
\leq\frac{C\log(n)}n,
\]
which is Item~(3) above. Finally, applying Lemma~\ref{lem4} to
Theorem~\ref{thm8} yields the claim
related to Item~(1) above
so that Theorem~\ref{thm7} is proved.
\end{pf*}

The remainder of this section is devoted to the proof of Theorem~\ref{thm8};
recall $W_{n,i}$ is the total degree of vertex $i$ in the preferential
attachment graph on $n$ vertices. Since we want to apply our geometric
approximation framework using the equilibrium distribution, we will use
Proposition~\ref{12} and so we first construct a variable having the size-bias
distribution of $W_{n,i}-1$. To facilitate this construction we need some
auxiliary variables.

For $j\geq i$, let $X_{j,i}$ be the indicator variable of the event
that vertex
$j$
has an outgoing edge connected to vertex $i$ in $G_j$
so that we can denote $W_{j,i}=1+\sum_{k=i}^j X_{k,i}$.
In this notation, for $1\leq i<j\leq n$,
%
\[
\IP(X_{j,i}=1\mid G_{j-1}) = \frac{W_{j-1,i}}{2j-1},
\]
and for $1\leq i \leq n$,
%
\[
\IP(X_{i,i}=1\mid G_{i-1}) = \frac{1}{2i-1}.
\]

The following well-known result will allow us to use this decomposition to
size-bias $W_{n,i}-1$; see, for example, Proposition~2.2 of \cite
{Chen2011} and the
discussion thereafter.

\begin{proposition}\label{prop2}
Let $X_1, \ldots, X_n$ be zero-one random variables with $\IP(X_i=1)=p_i$.
For each $i=1,\ldots, n$, let $(X_j^{(i)})_{j\not=i}$ have the
distribution of
$(X_j)_{j\not=i}$ conditional on $X_i=1$.
If $X=\sum_{i=1}^n X_i$, $\mu=\IE[X]$
and $K$ is chosen independently of the variables above with $\IP
(K=k)=p_k/\mu$,
then $X^s=\sum_{j\not=K} X_j^{(K)} +1$ has the size-bias distribution
of $X$.
\end{proposition}

Roughly, Proposition~\ref{prop2} implies that in order to size-bias $W_{n,i}-1$,
we choose an indicator
$X_{K,i}$ where, for $k=i,\ldots, n$, $\IP(K=k)$ is proportional to
$\IP(X_{n,k}=1)$ (and zero for other values),
then attach vertex
$K$ to vertex $i$, and sample the remaining edges conditionally on this event.
Note that given $K=k$, in the graphs $G_l$, $1\leq l < i$ and
$k<l\leq n$, this conditioning does not change the original rule for generating
the preferential attachment graph, given $G_{l-1}$.
The following lemma implies the remarkable fact that in order to
generate the
graphs $G_l$ for $i \leq l < k$
conditional on $X_{k,i}=1$ and $G_{l-1}$,
we attach edges following the same rule as preferential attachment, but include
the edge from vertex $k$
to vertex $i$ in the degree count.

\begin{lemma}\label{lem5}
Retaining the notation and definitions above,
for $i \leq j<k$, we have
%
\[
\IP(X_{j,i}=1\mid X_{k,i}=1, G_{j-1}) = \frac{1+W_{j-1,i}}{2j},
\]
where we define $W_{i-1,i}=1$.
\end{lemma}

\begin{pf}
By Bayes's rule, we have
%
\begin{equation} \label{29}
\IP(X_{j,i}=1\mid X_{k,i}=1, G_{j-1})
= \frac{\IP(X_{j,i}=1\mid G_{j-1}) \IP(X_{k,i}=1 \mid X_{j,i}=1,
G_{j-1})}{\IP(X_{k,i}=1\mid G_{j-1})},
\end{equation}
and we will calculate the three probabilities appearing in \eq{29}.
First, for $i\leq j$, we have
\[
\IP(X_{j,i}=1\mid G_{j-1}) = \frac{W_{j-1,i}}{2j-1},
\]
which implies
%
\[
\IP(X_{k,i}=1\mid G_{j-1}) = \frac{\IE[W_{k-1,i} \mid G_{j-1}]}{2k-1}
\]
and
%
\[
\IP(X_{k,i}=1\mid X_{j,i}=1, G_{j-1})
= \frac{\IE[W_{k-1,i} \mid X_{j,i}=1, G_{j-1}]}{2k-1}.
\]

Now, to compute the conditional expectations appearing above,
note first that
%
\[
\IE(W_{k,i}\mid G_{k-1}) =W_{k-1,i} + \frac{W_{k-1,i}}{2k-1}
= \biggl( \frac{2k}{2k-1} \biggr) W_{k-1,i},\vadjust{\goodbreak}
\]
and thus
%
\[
\IE(W_{k,i}\mid G_{k-2}) =\biggl( \frac{2(k-1)}{2(k-1)-1} \biggr)
\biggl( \frac{2k}{2k-1} \biggr) W_{k-2,i}.
\]
Iterating, we find that, for $i,s<k$,
%
\begin{equation}
\label{30}
\IE(W_{k,i}\mid G_{k-s})
=\prod_{t=1}^{s}\biggl( \frac{2(k-t+1)}{2(k-t+1)-1}
\biggr) W_{k-s,i}.
\end{equation}
By setting $j-1=k-s$ and then replacing $k-1$ by $k$ in \eq{30}, we obtain
%
\[
\IE(W_{k-1,i}\mid G_{j-1})
=\prod_{t=1}^{k-j}\biggl( \frac{2(k-t)}{2(k-t)-1} \biggr)
W_{j-1,i},
\]
which also implies
%
\[
\IE(W_{k-1,i}\mid X_{j,i}=1, G_{j-1})
=\prod_{t=1}^{k-j-1}\biggl( \frac{2(k-t)}{2(k-t)-1}
\biggr) (1+W_{j-1,i}).
\]
Substituting these expressions appropriately into \eq{29} proves the lemma.
\end{pf}

The previous lemma suggests the following (embellished) construction of
$(W_{n,i} | X_{k}=1)$ for any $1\leq i\leq k \leq n$. Here and below we will
denote quantities related to this construction with a superscript $k$.
First we
generate $G_{i-1}^k$, a graph with $i-1$ vertices, according to the usual
preferential attachment model. At this point, if $i\not=k$, vertex $i$
and $k$
are added to the graph, along with a vertex labeled $i'$ with an edge
to it
emanating from vertex $k$. Given $G_{i-1}^k$ and these additional
vertices and
edges, we generate $G_i^k$ by connecting vertex $i$ to a vertex $j$ randomly
chosen from the vertices $1, \ldots, i, i'$ proportional to their
degree, where
$i$ has degree one (from the out-edge), and $i'$ has degree one (from
the in-edge
emanating from vertex $k$). If $i=k$, we attach $i$ to $i'$ and denote the
resulting graph by $G_i^i$. For $i < j <k$, we generate the graphs $G_{j}^k$
recursively from $G_{j-1}^k$ by connecting vertex $j$ to a vertex $l$, randomly
chosen from the vertices $1, \ldots, j, i'$ proportional to their
degree, where
$j$ has degree one (from the out-edge). Note that none of the vertices $1,
\ldots, k-1$ can connect to vertex $k$. We now define
$G_k^k=G_{k-1}^k$, and for
$j=k+1, \ldots, n$, we generate $G_j$ from $G_{j-1}$, according to preferential
attachment among the vertices $1, \ldots, j, i'$.

\begin{lemma} \label{lem6}
Let $1\leq i \leq k \leq n$, and retain the notation and definitions above.
\begin{enumerate}[(3)]
\item[(1)] $\law(W^{k}_{n,i}+W^{k}_{n,i'}) = \law(W_{n,i} | X_{k}=1).$\vspace*{1pt}
\item[(2)] For fixed $i$, if $K$ is a random variable such that
%
\[
\IP(K=k) = \frac{\IE X_{k,i}}{\IE W_{n,i}-1},\qquad k\geq i,
\]
then $W^{K}_{n,i}+W^{K}_{n,i'}-1$ has the size-bias distribution of $W_{n,i}-1$.
\item[(3)] Conditional on the event $\{W^{k}_{n,i}+W^{k}_{n,i'}=m+1\}$,
the variable $W^{k}_{n,i}$ is uniformly distributed on the integers $1, 2,
\ldots, m$.
\item[(4)] $W^{K}_{n,i}-1$ has the discrete equilibrium distribution of
$W_{n,i}-1$.
\end{enumerate}
\end{lemma}
\begin{pf}
Items (1) and (2) follow from Proposition~\ref{prop2} and Lemma~\ref{lem5}. Viewing
$(W^{k}_{n,i},W^{k}_{n,i'})$ as the number of balls of two colors in a
Polya urn
model started with one ball of each color, Item~(3) follows from
induction on $m$
and Item (4) follows from Proposition~\ref{12}.
\end{pf}

\begin{pf*}{Proof of Theorem~\protect\ref{thm8}}
We will apply Theorem~\ref{thm2} to $\law(W_{n,i}-1)$, so that we must
find a coupling of a variable with this distribution to that of a variable
having its discrete equilibrium distribution. For each fixed $k=i,
\ldots, n$
we will construct ${(X^{k}_{j,i}, \widetilde{X}^{k}_{j,i})}_{j\geq
i}$ so
${(X^{k}_{j,i})}_{j\geq i}$ and ${(\widetilde
{X}^{k}_{j,i})}_{j\geq i}$
are distributed as the indicators of the events
vertex $j$ connects to vertex $i$ in $G_n^k$ and~$G_n$, respectively.
We will use the notation
%
\[
W_{j,i}^k=\sum_{m=i}^j X_{j,i}^k
\quad\mbox{and}\quad
\widetilde{W}_{j,i}^k=\sum_{m=i}^j \widetilde{X}_{j,i}^k,
\]
which will be distributed as the degree of vertex $i$ in the appropriate
graphs.

The constructions for $k=i$ and $k>i$ differ, so assume here that $k>i$.
Let $U_{j,i}^k$ be independent uniform $(0,1)$ random variables,
and first define
%
\[
X_{i,i}^k=\I{[U_{i,i}^k<1/2i]}
\quad\mbox{and}\quad
\widetilde{X}_{i,i}^k=\I{[U_{i,i}^k<1/(2i-1)]}.
\]
Now, for $i<j<k$, and assuming that ${(W^k_{j-1,i},
\widetilde{W}^k_{j-1,i})}$ is given, we
define
%
\begin{equation} \label{31}
X_{j,i}^k=\I{\biggl[U_{j,i}^k<\frac{W^k_{j-1,i}}{2j}\biggr]}
\quad\mbox{and}\quad
\widetilde{X}_{j,i}^k=\I{\biggl[U_{j,i}^k<
\frac{\widetilde{W}^k_{j-1,i}}{2j-1}\biggr]}.
\end{equation}
For $j=k$ we set $X_{k,i}^k=0$ and $\widetilde{X}_{k,i}^k$ as in \eq
{31} with
$j=k$,
and, for $j>k$, we define
%
\[
X_{j,i}^k=\I{\biggl[U_{j,i}^k<\frac{W^k_{j-1,i}}{2j-1}\biggr]}
\quad\mbox{and}\quad
\widetilde{X}_{j,i}^k=\I{\biggl[U_{j,i}^k
<\frac{\widetilde{W}^k_{j-1,i}}{2j-1}\biggr]}.
\]
Thus we have recursively defined the variables ${(X^{k}_{j,i},
\widetilde{X}^{k}_{j,i})}$ and it is clear they are distributed as claimed
with ${(W_{j,i}^k, \widetilde{W}_{j,i}^k)}$ distributed as the required
degree counts. Note also that
$\widetilde{W}_{j,i}^k\geq W_{j,i}^k$ and $\widetilde{X}_{j,i}^k\geq X_{j,i}^k$.
We also define the events
%
\[
A_{j,i}^k := \bigl\{\min\{i\leq l\leq n\dvt
X_{l,i}^k\not=\widetilde{X}^k_{l,i}\}=j\bigr\}.
\]
Using that $W_{j-1,i}^k=\widetilde{W}^k_{j-1,i}$ under $A_{j,i}^k$
(which also
implies $A_{j,i}^k=\varnothing$ for $j>k$) we have
\begin{eqnarray*}\label{32}
\IP{(\widetilde{W}^{k}_{n,i} \neq{W}^{k}_{n,i})}
& = & \IP{\Biggl( \bigcup_{j=i}^n A_{j,i}^k\Biggr) } \\
&= &\IE\widetilde{X}^k_{k,i}+ \sum_{j=i}^k
\IP\biggl( A_{j,i}^k\cap\biggl\{\frac{W_{j-1,i}^k}{2j}
< U_{j,i}^k < \frac{\widetilde{W}_{j-1,i}^k}{2j-1}\biggr\}\biggr) \\
&\leq&\IE\widetilde{X}^k_{k,i}+ \sum_{j=i}^{n}
\IP\biggl( \frac{\widetilde{W}_{j-1,i}^k}{2j} < U_{j,i}^k
< \frac{\widetilde{W}_{j-1,i}^k}{2j-1}\biggr),
\end{eqnarray*}
where\vspace*{1pt} we write $W_{i-1,i}^k:=\widetilde{W}_{i-1,i}^k:=1$.
Finally, starting from \eq{32}
and using
the computations in the proof of Lemma~\ref{lem5} and
the estimates in Lemma~\ref{lem3}, we find
\begin{eqnarray*}
\IP{(\widetilde{W}^{k}_{n,i} \neq{W}^{k}_{n,i})}
&\leq&\frac{\IE W_{k-1,i}}{2k-1} +\sum_{j=i}^{n} \IE W_{j-1,i}
\biggl( \frac{1}{2j-1}-\frac{1}{2j} \biggr) \\
&\leq& C{\Biggl[\sqrt{\frac{k}{i}}\frac{1}{k}+
\sqrt{\frac{k}{i^3}}\frac{1}{k}
+\sum_{j\geq i} \Biggl(\sqrt{\frac{j}{i}}\frac{1}{j^{2}}+
\sqrt{\frac{j}{i^{3}}}\frac{1} {j^{2}}\Biggr)\Biggr]} \leq C/i.
\end{eqnarray*}
If $k=i$, it is clear from the construction preceding Lemma~\ref{lem6}
that an
easy coupling, similar to that above,
will yield $\IP(\widetilde{W}^{i}_{n,i} \neq{W}^{i}_{n,i}) <C/i$.
Since these bounds do not depend on $k$, we also have
\begin{equation}
\label{33}
\IP{(\widetilde{W}^{K}_{n,i}-1 \neq{W}^{K}_{n,i}-1)} \leq C/i,
\end{equation}
and the result now follows from Lemma~\ref{lem6}, \eq{33} and~\eq{11} of
Remark~\ref{rem3}.
\end{pf*}

\section*{Acknowledgements}
Erol A. Pek\"oz and
Adrian R\"ollin were supported in part by NUS Grant R-155-000-098-133.
The authors would like to thank the anonymous referees for their helpful
comments.
%

\printhistory

\end{document}